\documentclass[submission,nomarks]{dmtcs-episciences}
\usepackage[round]{natbib}
\usepackage[usenames,x11names]{xcolor}
\usepackage{amssymb,amsmath,amsthm,amsfonts,amscd}
\usepackage[amsstyle]{cases}
\usepackage{graphicx}
\definecolor{purp}{rgb}{0.5, 0.0, 0.5}
\newcommand{\purple}{\color{purp}}
\newcommand{\teal}{\color{teal}}
\usepackage{latexsym,accents,fdsymbol,stmaryrd}
\usepackage[pdf]{pstricks}
\usepackage{mathtools,epsf,pst-node,pst-tree,auto-pst-pdf,pst-plot}
\newtheorem{theorem}{Theorem}
\newtheorem{lemma}[theorem]{Lemma}
\newtheorem{corollary}[theorem]{Corollary}
\newcommand{\AH}[4]{\smash{{#1}\mathop{\longleadsto}\limits^{\raisebox{-1pt}{$\scriptscriptstyle #3$}}_{\smash{\raisebox{1mm}{$\scriptscriptstyle #4$}}}{#2}}}
\newcommand{\A}[3]{\AH{#1}{#2}{#3}{}}

\newcommand{\brk}[1]{\llbracket #1 \rrbracket}
\newcommand{\br}[2]{\brk{#1\pm #2}}
\newcommand{\D}[3]{\FFF{#1}{k}{#2}{#3}{h}}    
\newcommand{\GH}[4]{\AH{#1}{[#2]}{#3}{#4}}
\newcommand{\G}[3]{\GH{#1}{#2}{#3}{}}
\newcommand{\FF}[4]{\FFF{#1}{#2}{#3}{#4}{}} 
\newcommand{\FFF}[5]{\AH{#1}{\br{#2}{#3}}{#4}{#5}}

\newcommand{\NE}{\ensuremath\nearrow}
\newcommand{\SE}{\ensuremath\searrow}

\newcommand{\dotdot}{\mathbin :} 
\usepackage[version-1-compatibility]{siunitx}
\newcolumntype{d}{S[tabformat=4.0,tabautofit]}   
\newcommand{\floor}[1]{{\left\lfloor #1\right\rfloor}}
\newcommand{\ceil}[1]{{\left\lceil #1\right\rceil}}
\renewcommand{\vec}[1]{\accentset{\longleftarrow}{#1}}
\newcommand{\seqnum}[1]{\href{http://oeis.org/#1}{{#1}}}
\usepackage{epigraph}
\title[Between Broadway and the Hudson: A Bijection of Corridor Paths]{Between Broadway and the Hudson:\\A Bijection of Corridor Paths}
\author{Nachum Dershowitz}
\affiliation{School of Computer Science, Tel Aviv University, Ramat Aviv, Israel}
\date{December 2020}       
\begin{document}
\makeatletter\gdef\@firsthead{}\gdef\@copyright{}\makeatother  
\maketitle

\epigraph{\small Canal street, running across Broadway to the Hudson,
near the centre of the city,\\
is a spacious street, principally occupied by retail stores\dots.\\
The streets are generally well paved,
with good side walks,\\
lighted at night with lamps,
and some of them supplied with gas lights.}{---\textit{The Treasury of Knowledge, and Library of Reference} (1834)}

\begin{abstract}
We present a substantial generalization of the equinumeracy of Grand Dyck paths and Dyck-path prefixes, constrained within a band.
The number of constrained paths starting at level $i$ and ending in a window of size $2j+2$ is
equal to the number starting at level $j$ and ending in a window of size $2i+2$ centered around the same point.
A new encoding of lattice paths provides a bijective proof.
\end{abstract}

\section{Introduction}

\begin{figure}[t!]
\begin{center}
\includegraphics[width=0.75\linewidth,trim= 0 100 0 80,clip]{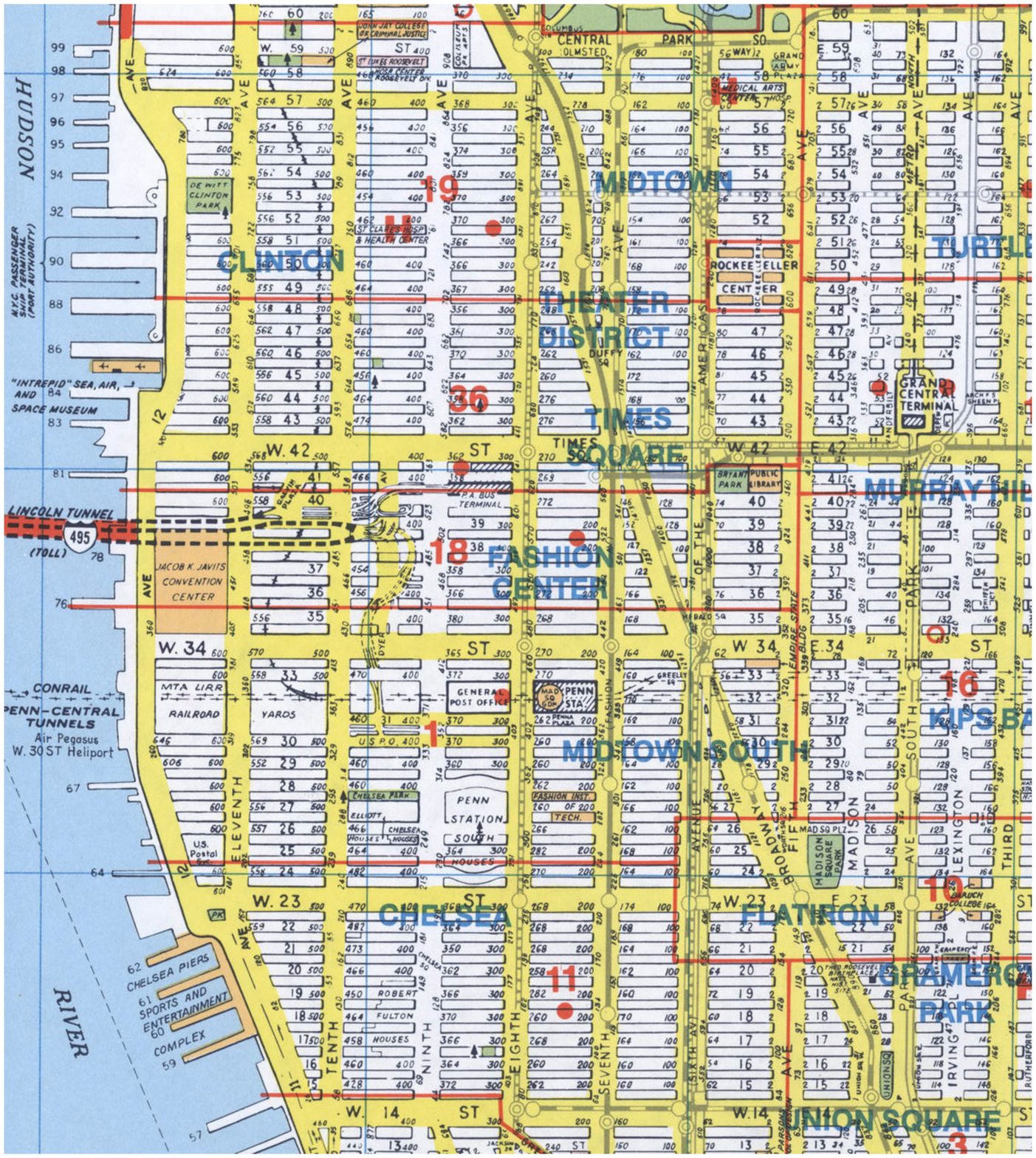}
\end{center}
\caption{Manhattan neighborhoods, with east-west streets and north-south avenues, bounded by Broadway on the east and the Hudson River on the west, with Union Square serving as origin.
(Image \raisebox{2pt}{\tiny\copyright}\;Hagstrom Map Company, Inc., in the public domain at \url{https://www.maps-of-the-usa.com/usa/new-york/new-york/large-detailed-road-map-of-south-manhattan-nyc}.)}\label{fig:bdwy}
\end{figure}

It has long been known that the Grand Dyck lattice paths (corresponding to ballot counting ending in a tie) and Dyck-path prefixes (when one candidate remains in the lead throughout the counting) are in bijection \citep[p.\@ 96]{Feller}.
More recently it has been shown that this is also true when steps are restricted to a band of restricted height \cite[]{Cigler,GP},
sometimes called \emph{corridor paths} \citep[e.g.\@][]{AK}.
Imagine walking in Manhattan, sticking west of Broadway (Figure~\ref{fig:bdwy}).
We set out to generalize this relationship between sets of paths to all possible starting points and ending ranges.

Let 
$\AH{i}{\ell}nh$, or just $\A{i}{\ell}n$ (fixing $h$),
denote the number of monotonic lattice paths from $\langle 0,i\rangle$ to $\langle n,\ell\rangle$
with $n$ steps that stay within (but may touch) the boundaries $y=0$ and $y=h$, for some given (maximum) height $h$.\footnote{\emph{Height} here is the maximum length of a unidirectional path (just NE or just SE).
Some  prefer to say that the \emph{width} of the corridor is $h+1$, since
$h+1$ ordinate values are allowed.}
Let $H=[0\dotdot h]=\{0,1,\ldots,h\}$ be the ordinate bounds within which steps are permissible.
Steps are diagonal, NE (northeast, $\nearrow$, a ``rise''), taking $\langle x,y\rangle\mapsto\langle x+1,y + 1\rangle$, and SE (southeast, $\searrow$, a ``fall''), taking $\langle x,y\rangle\mapsto\langle x+1,y - 1\rangle$, both with the proviso that the new ordinate position $y\pm 1\in H$, as the case may be.  
It is easy to see that the only nonzero counts of  $n$-step paths starting at level $i$ and ending at $\ell$ occur when
$n +i \equiv \ell \pmod 2$.
See Figure~\ref{fig:path} for a sample path in $\AH{1}{3}{12}4$ (using the same notation  for the set of paths as for its cardinality).

\begin{figure}[t!]
\centering
\begin{pspicture}(0,-0.5)(12,5)
\psgrid[subgriddiv=1,griddots=10](0,0)(0,0)(12,4)
\uput{10pt}[150](-0.1,3.8){$h=$}
\uput{10pt}[150](-0.1,1.8){$k=$}
\uput{10pt}[150](-0.1,0.8){$i=$}
\uput{5pt}[0](12,3){$\ell$}
\uput{15pt}[0](12,2.5){$\left.\rule{0pt}{50pt}\right\}=J$}
\psset{linecolor=red,linewidth=2pt,linestyle=solid}
\pcline{-}(0,0)(12,0)
\pcline{-}(0,4)(12,4)
\pcline[linestyle=dashed,,linewidth=1pt,linecolor=red]{-}(0,2)(12,2)
\psset{linecolor=blue,arrowsize=6pt}
\pcline[linestyle=solid]{->}(0,1)(1,2)
\pcline[linestyle=solid]{->}(1,2)(2,3)
\pcline[linestyle=solid]{->}(2,3)(3,2)
\pcline[linestyle=solid]{->}(3,2)(4,3)
\pcline[linestyle=solid]{->}(4,3)(5,4)
\pcline[linestyle=solid]{->}(5,4)(6,3)
\pcline[linestyle=solid]{->}(6,3)(7,2)
\pcline[linestyle=solid]{->}(7,2)(8,1)
\pcline[linestyle=solid]{->}(8,1)(9,0)
\pcline[linestyle=solid]{->}(9,0)(10,1)
\pcline[linestyle=solid]{->}(10,1)(11,2)
\pcline[linestyle=solid]{->}(11,2)(12,3)
\end{pspicture}
\medskip
\caption{A diagonal path 
(counted by)
 {\blue $\A{1}{3}{12}$}, 
which goes from $i=1$ to $\ell=3$ in a dozen steps,  consisting of 7 NE-steps and 5 SE-steps, with bound $h=4$.
The target region $J$ is $\br{2}{1}=[1\dotdot 4]$;
its center  is $k=2$ (the dashed {\red red} line), and $j=1$ determines its width.
There is one \emph{feasible} endpoint at or below $k$ and one above.}\label{fig:path}
\end{figure}
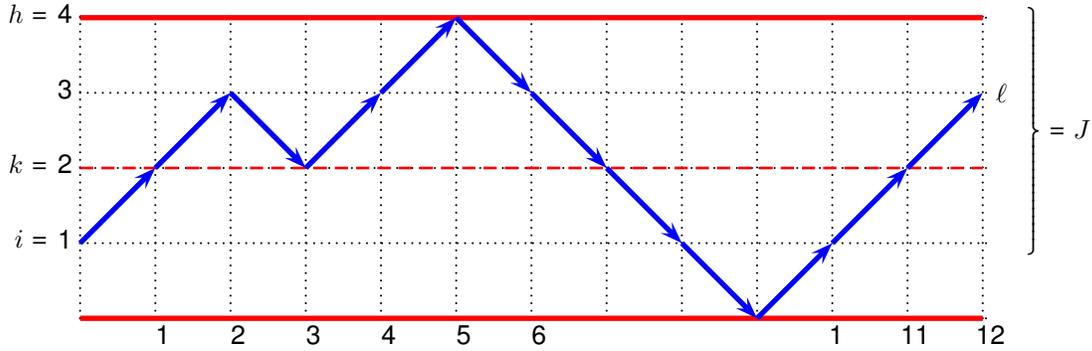

\begin{figure}[t!]
\centering
\begin{pspicture}(-0.2,-0.2)(10,1.6)
\psset{xunit=2.5cm,yunit=2.5cm,linecolor=darkgray,linestyle=solid,linewidth=2pt}
\psgrid[subgriddiv=1,griddots=20,gridlabels=0pt,gridwidth=2pt](0,0)(0,0)(4,0)
\rput(0,-0.1){\scriptsize 0}
\rput(1,-0.1){\scriptsize 1}
\rput(2,-0.1){\scriptsize 2}
\rput(3,-0.1){\scriptsize 3}
\rput(4,-0.1){\scriptsize 4}
\psline[]{o->}(1,0)(2,0)
\psline{->}(2,0)(3,0)
\psline[linecolor=teal]{>->}(3,0.05)(2,0.1)
\psline{>->}(2,0.15)(3,0.2)
\psline{->}(3,0.2)(4,0.2)
\psline[linecolor=teal]{>->}(4,0.25)(3,0.3)
\psline[linecolor=teal]{->}(3,0.3)(2,0.3)
\psline[linecolor=teal]{->}(2,0.3)(1,0.3)
\psline[linecolor=teal]{->}(1,0.3)(0,0.3)
\psline{>->}(0,0.35)(1,0.4)
\psline{->}(1,0.4)(2,0.4)
\psline{-**}(2,0.4)(3,0.4)
\pcline[linestyle=dashed,linecolor=Green3,linewidth=2pt]{-}(2.5,-0.2)(2.5,0.6)
\end{pspicture}
\medskip
\caption{A right-left version of the constrained path $\AH{1}{3}{12}4$, consisting of 7 right ($\mathsf{\darkgray +}$) steps (colored {\darkgray gray}) and 5 left ($\mathsf{\teal -}$) steps ({\teal teal}) along a 5-vertex point graph $P_5$ (labeled 0,1,2,3,4), starting at vertex $i=1$.
The path in this representation is $\mathsf{\darkgray +}\mathsf{\darkgray +}\mathsf{\teal -}\mathsf{\darkgray +}\mathsf{\darkgray +}\mathsf{\teal -}\mathsf{\teal -}\mathsf{\teal -}\mathsf{\teal -}\mathsf{\darkgray +}\mathsf{\darkgray +}\mathsf{\darkgray +}$, based at $1$.
It is an accordion fold of the blue path in Figure~\ref{fig:path}.
The {\color{Green4} green} vertical line serves as a ``center of attraction.'' 
(See Section~\ref{sec:bi} and Figure~\ref{fig:ta}.)}\label{fig:line}
\end{figure}

The basic recurrence is
\begin{align*}
\A i\ell n &= \begin{cases}
0 &\mbox{if } i\notin H \mbox{ or } \ell\notin H\\
[i=\ell] &\mbox{if } n=0\\
(\A{i}{\ell-1}{n-1}) +( \A{i}{\ell+1}{n-1})  &\mbox{otherwise}
\end{cases}
\end{align*}
where the bracketed condition $[i=\ell]$ is Iverson's notation for a characteristic function (1 when true; 0 when false),
and the conditions are taken in order.

The ends of the paths we are interested in  fall within a  range, $J$, not just a single point $\ell$.
For example, the window $J=[5\dotdot 10]$
has $6$ possible landing spots, but only half of them are feasible, depending on
whether $n+i$ is odd or even.
Only those $\ell\in J$ with the same parity as $n+i$ are relevant.
Our goal is to count
\begin{align*}
\A{i}{J}n &=\sum_{\ell\in J} \A{i}{\ell}n ~=~ \sum_{\clap{$\substack{\ell\in J\\\ell\equiv n+i \!\!\!\!\pmod 2}$}} \A{i}{\ell}n 
\end{align*}
the number of paths constrained to any corridor $H=[0\dotdot h]$ and ending at any (feasible) ordinate in the  window $J$.

These constrained lattice paths are equivalent to walks along a path graph, forward and backward.
See Figure~\ref{fig:line}.
When $i=0$ (at the bottom) and $J=H$ (anywhere), walks for $h=0,1,2,3,4,5,6,7,8,9,10$ are enumerated at
\seqnum{A000007} (constant 0 after initial 1), 
\seqnum{A000012} (constant\nolinebreak{} 1), 
\seqnum{A016116} (\smash{$2^\floor{n/2}$}), 
\seqnum{A000045} (Fibonacci), 
\seqnum{A038754} ($\{1,2\}3^n$), 
\seqnum{A028495}, 
\seqnum{A030436}, 
\seqnum{A061551},
\seqnum{A178381}, 
\seqnum{A336675},
\seqnum{A336678},
respectively,
in Neil Sloane's \textit{Encyclopedia of Integer Sequences (OEIS)} \cite[]{OEIS}.\footnote{Compiled largely by Jonathon Bryant \cite[]{AK}. 
The sequences for $h=9,10$ were added recently.}

Such paths in a path graph having $h$ edges can also be viewed as prefixes of Dyck paths of bounded height $h$,
since they start at the bottom but may end anywhere above or on the bottom line.
Their number is known to be equal to that of Grand Dyck paths, of the same length, which start in the middle of the band, may go above or below that line~-- as long as they stay within bounds,
and which we allow to end up either in the middle or just above~\cite[]{Cigler}.%
\footnote{Usually, Grand Dyck paths are defined to be of even length and to end up back on the starting line
(as in a tie vote).
To be more inclusive, we allow odd-length Grand Dyck paths that terminate one line above --
as in \citet{Elizalde2015BijectionsFP}, for instance~-- adopting the same moniker in the odd case, too.
Accordingly, we can say that the number of  Grand Dyck paths of length $2n$ (``semilength'' $n$)  and even height $h$ (with the up-down symmetry of the corridor) is always twice that for $2n-1$.
See the middle case of Table~\ref{tab}.
The term ``Grand Dyck'' is used in \citet{Shapiro}, for example;
these lattice paths are also referred to 
as ``two-sided'' or ``bilateral'' paths \citep[e.g.\@][]{LiptonZ77} on account of their shape,
as ``binomial'' or ``central binomial'' paths  \citep[e.g.\@][]{Pergola} on account of their number, and
as  ``{free}'' Dyck paths \citep[e.g.\@][]{CPQS}.
These paths are mentioned in \cite[p.\@ 81]{Comtet74}.
They are classified as ``bridges'' in \citet{CM,BF}.}
So, each of the above sequences also counts constrained Grand Dyck paths.
This explains the inclusion also of even- and odd-indexed partitions of the corridor sequences in the \textit{OEIS}:
they are enumerations of (combinatorial objects in bijection with) the Grand Dyck paths ending on the diagonal (in the middle of the corridor) and of those ending one row above.
For example, sequence
\seqnum{A336678} for $h=10$ is the alternation of sequences \seqnum{A087944} and \seqnum{A087946}.

More generally, walks can start anywhere in $H$ ($0\leq i\leq h$), 
with the ordinal positions along the route always staying within $H$.
Table \ref{tab} (at the end) lists values for the number of paths through a corridor of height $h=4$, with one subtable for each starting point ($i=0,1,2,3,4$);
Table \ref{tab:5} exhibits $h=5$.
These may be viewed as constrained versions of Pascal's triangle, with each entry the sum of two prior entries
\citep[cf.\@][]{AK}.

In addition to pointing to the formula (Theorem~\ref{thm:cnt}) enumerating these more general sets of corridor paths ending in an arbitrary window,
we explore a beautiful symmetry between such sets of paths,
those starting at level $i$ and ending in a window of size $2j+2$ and those
starting at level $j$ and ending in a window of size $2i+2$ centered around the same point
(Theorem~\ref{thm:sym}).
Three proofs are then provided for this symmetry: by induction (Section~\ref{sec:ind}), by counting (Section~\ref{sec:comb}), and by bijection (Section~\ref{sec:bi}).
The final section describes prior work leading up to these results.
For a brief history of lattice-path enumerations, see \cite[]{history}.

\section{Main Results}

We use the notation $\br{k}{j}$
as shorthand for a range $[k-j\dotdot k+j+1]$,
which we make of even size, viz.\@ $2j+2$, by stretching the upper end one spot, to include $k+j+1$.
Thus,
the window $\br{k}{j}$ covers $j+1$ feasible endpoints~-- the odd ones or the even ones, as the case may be~-- centered about $k$.

Our main result is the following intriguing equivalence:
\begin{theorem}\label{thm:sym}
For all $n,h\in\naturals$, $k\in [0\dotdot h]$, $i,j\in[0\dotdot \min\{k+1,h-k\}]$: 
\begin{align}\label{eq:sym}
\FFF{i}{k}{j}{n}h & ~~=~~ \FFF{j}{k}{i}{n}h
\end{align}
\end{theorem}
\noindent
For example, $\FFF{2}{2}{1}{9}4  = 162 = \FFF{1}{2}{2}{9}4$;
see Table \ref{tab}.
The bounds on $i$ and $j$ ensure that the starting points are in $H=[0\dotdot h]$ and that the target windows $\br{k}j$ and $\br{k}i$ do not extend beyond one row above or below the corridor $H$.
(The reason for allowing the window to slightly overextend the corridor is that~-- when $h$ is even~-- there are an odd number of possible points in the full height of the corridor, all of which one may wish to include, whereas target windows necessarily cover an even number of landing points.)

Were $i$ or $j$ too big, $k\pm i$  or $k\pm j$ could extend too far beyond $H$, and the equality would not hold,
as is the case for  $\FFF{2}{3}{1}{9}4  = 81 \neq \FFF{1}{3}{2}{9}4 = 121$.
When  $k\leq h/2$,  the theorem  holds  as long as $i,j\leq k$.

The largest $i$ and $j$ can be is $i=\floor{h/2}$ and $j=\ceil{h/2}$ (or vice versa), 
which gives
\begin{align}\label{case}
 \A{\floor{h/2}}{H}{n} & ~~=~~ \A{\ceil{h/2}}{H}{n}
\end{align}
and comes as no surprise.

This theorem also holds for the degenerate case $k=-1$
since the constraints impose $i=j=0$, in which case the equivalence is true trivially,
and the paths are Dyck paths of bounded height.

By up-down symmetry:
\begin{lemma}\label{lem}
For all $n,h\in\naturals$, $i,j,k\in [0\dotdot h]$,
\begin{align}
\FFF{i}{k}{j}{n}h & ~~=~~ \FFF{(h-i)}{h-k-1}{j}{n}h
\end{align}
\end{lemma}

So, for instances when $i>h\div 2$ ($=\lfloor h /2\rfloor$, using the obelus for integer division), we can combine this lemma with our theorem to obtain:
\begin{corollary}\label{cor:2}
The equivalence
\begin{align*}
 \FFF{i}{k}{j}{n}h & ~~=~~ \FFF{j}{h-k-1}{h-i}{n}h
\end{align*}
holds for all $n,h\in\naturals$,  $k\in [0 \dotdot h]$, 
   $i\in [\max\{k,h-k-1\}\dotdot h]$,
 $j\in[0\dotdot \min\{k+1,h-k\}]$.
\end{corollary}

If no upper boundary on paths is imposed (effectively when $h\geq n$), then we have the following:
\begin{corollary}
For all $n,k\in\naturals$, $i,j\in[0\dotdot k+1]$: 
\begin{align}\label{eq:inf}
\FFF{i}{k}{j}{n}\infty & ~~=~~ \FFF{j}{k}{i}{n}\infty
\end{align}
\end{corollary}
\noindent
The equivalence of unrestricted Grand Dyck and Dyck prefix paths is just one special case
($i=k \geq n/ 2$, $j=0$).

In terms of ballots, this means that
the likelihood of candidate A never being behind B during the counting by more than $j$ votes and then winning by $\ell\mathbin{..} \ell+2i+1$ votes (losing, when negative) is the same as 
the likelihood of  A winning (or losing) by $\ell-j\mathbin{..} \ell+j+1$ and never being behind B during the counting by more than $i$ ($\ell+1\geq j-i,-j$).
When $j=0$, this translates into the following:
\begin{corollary}
The likelihood of candidate A never being behind candidate B by more than $i$ votes during the counting of ballots and winning by exactly $\ell$ votes 
 is equal to 
the likelihood of  A never being behind at all and winning by $\ell$ to $\ell+2i$ votes (for all feasible values of $\ell$).
\end{corollary}
\noindent
Feasible means that $\ell\geq 0$ and has the right parity (same as the number of steps) for an outcome.

Lastly, a closed-form formula for the paths of interest is as follows:
\begin{theorem}\label{thm:cnt}
The number of corridor paths $\D{i}{j}n$, for all $n,h,i,j,k\in\naturals$, $i\leq h$, is
\begin{align*}
\sum_{z=v}^u
\sum_{\substack{s=0\\0\leq 2s+p \leq h}}^{j}
\Bigg[ \binom{n}{r+z(h+2)+s}  - \binom{n}{r+z(h+2)+i+s+1} \Bigg]
\end{align*}
where $q=n+k-i-j$, $r=\ceil{q/2}$, $p=k-j+q\bmod 2$,
$v=-\floor{(r+i+j+1)/(h+2)}$ and $u=\floor{(n-r)/(h+2)}$.
\end{theorem}

\section{Inductive Proof}\label{sec:ind}

One can prove Theorem \ref{thm:sym}, namely that
\begin{align*}
\FFF{i}{k}{j}{n}h & = \FFF{j}{k}{i}{n}h
\end{align*}
by induction on the number of steps $n$,  with height $h$ fixed throughout.

Recall that the bounds on  $i$ and $j$ are
\begin{align}
\label{4} i,j &\geq 0\\
\label{5} i,j &\leq k+1\\
\label{6} i,j &\leq h-k
 \end{align}
The cases where either is out of bounds 
are excluded from the theorem.

For $n=0$, the starting and ending points must be the same.
The two boundary conditions, viz.
\begin{align*}
\A i{\br k j} 0 &=[ k-j\leq i\leq k+j+1]\\
\A j{\br k i} 0 &= [ k-i\leq j\leq k+i+1]
\end{align*}
are equivalent
since we are given that $0\leq i,j\leq k+1$.

In the general case ($n>0$), we could argue inductively in the following fashion:
\begin{align*}
& \A{i}{\br k j}{n} \\
{}={}
& (\A{i-1}{\br k j}{n-1}) +(\A{i+1}{\br k j}{n-1})  & \mbox{basic recurrence} \\
{}={}
&  (\A{j}{\br k {i-1}}{n-1}) +(\A{j}{\br k {i+1}}{n-1}) & \mbox{induction} \\
{}={}
&  (\A{j}{\br {k-1} {i}}{n-1}) +(\A{j}{\br {k+1} {i}}{n-1}) & \mbox{definition} \\
{}={}
& \A{j}{\br k i}{n}                & \mbox{basic recurrence} 
\end{align*}
But this  only works if the two inductive cases also satisfy the theorem's constraints.

The problematic cases, when the inductive hypothesis cannot be applied, are three: 
\begin{itemize}
\item[(a)] $i=0$, since then $i-1<0$ violates (\ref{4}) for the left inductive case $\A{i-1}{\br k j}{n-1}$; 
\item[(b)] $i=k+1$, since then $i+1>k+1$ in violation of (\ref{5}) for the right inductive case $\A{i+1}{\br k j}{n-1}$; 
\item[(c)] $i=h-k$, violating (\ref{6}) for the right case.
\end{itemize}

Fortuitously,
the exact same argument may be applied in the opposite direction, with the r\^oles of $i$ and $j$ exchanged, to prove the identical equivalence:
\begin{align}\label{j}
\A{i}{\br k j}{n}         
{}={}
&  (\A{i}{\br k {j-1}}{n-1}) +(\A{i}{\br k {j+1}}{n-1})  \\
\nonumber
{}={}
& (\A{j-1}{\br k i}{n-1}) +(\A{j+1}{\br k i}{n-1})  \\
\nonumber
{}={}
& \A{j}{\br k i}{n} 
\end{align}
The  cases for which this version of the argument is problematic are  analogous but different:
\begin{itemize}
\item[(a')] $j=0$; 
\item[(b')]
 $j=k+1$;  or 
 \item[(c')] $j=h-k$.
 \end{itemize}

For the first exception (a), when  $i=0$, all is well with just one induction:
\begin{align*}
\A{0}{\br k j}{n}
= \A{1}{\br k j}{n-1} 
&= 
 \A{j}{\br k 1}{n-1}
= \A{j}{\br k 0}{n}
\end{align*}
In the extreme case that $k=h$, and the induction is invalid, it must also be that $j=0$,
and the equivalence holds immediately, sans induction.
By the same token, case (a') is also not an issue.

Furthermore, whenever $i=j$, the theorem holds trivially, so the two combined
cases (b,b'), when $i=j=k+1$, and (c,c'), when $i=j=h-k$, are fine, too.

So we only lack a proof for the following two combinations of the exceptions:
(b,c'), when  $i=k+1$, $j=h-k$, and (c,b'), when $j=k+1$ and $i=h-k$.
These are symmetric, so let's  delve just into the second.
Taking constraints (\ref{5},\ref{6}) into account, we find that $h=2k+1$, 
$i=k=\floor{h/2}$, and $j=\ceil{h/2}$.
So all we have to establish is the case
$\A{\floor{h/2}}{H}{n}  = \A{\ceil{h/2}}{H}{n}$, which we've already seen 
(\ref{case}).

\section{Combinatorial Proof}\label{sec:comb}

\begin{figure}[t!]
\centering
\begin{pspicture}(4,8)
\psgrid[subgriddiv=1,griddots=10](0,0)(0,0)(5,7)
\psset{linecolor=red,linewidth=2pt,linestyle=solid}
\psline[linestyle=solid]{-}(2,0)(5,3)
\psline[linestyle=solid]{-}(0,4)(3,7)
\psline[linestyle=dashed]{-}(-0.4,3.6)(0,4)
\psline[linestyle=dashed]{-}(1.6,-0.4)(2,0)
\psline[linestyle=dashed]{-}(5,3)(5.4,3.4)
\psline[linestyle=dashed]{-}(3,7)(3.4,7.4)
\psset{linecolor=blue,arrowsize=6pt}
\psline[linestyle=solid]{->}(0,0)(0,1)
\psline[linestyle=solid]{->}(0,1)(0,2)
\psline[linestyle=solid]{->}(0,2)(1,2)
\psline[linestyle=solid]{->}(1,2)(1,3)
\psline[linestyle=solid]{->}(1,3)(1,4)
\psline[linestyle=solid]{->}(1,4)(2,4)
\psline[linestyle=solid]{->}(2,4)(3,4)
\psline[linestyle=solid]{->}(3,4)(4,4)
\psline[linestyle=solid]{->}(4,4)(5,4)
\psline[linestyle=solid]{->}(5,4)(5,5)
\psline[linestyle=solid]{->}(5,5)(5,6)
\psline[linestyle=solid]{->}(5,6)(5,7)
\uput{10pt}[150](-0.0,3.77){$s{=}$}
\uput{10pt}[150](2.2,-0.5){$t{=}$}
\uput{10pt}[150](6.4,7.0){$\langle a , b\rangle$}
\end{pspicture}
\medskip
\caption{An orthogonal (``staircase'') path, starting at the origin and ending at $\langle a , b\rangle=\langle 5,7\rangle$, consisting of 7 N-steps and 5 E-steps, staying strictly within bounds $s=4$ (below {\red $y=x+4$}) and $t=2$  (above {\red $y=x-2$}).  This path and its constraints are analogues of those in Figure~\ref{fig:path};
see Section \ref{sec:comb}.}
\label{fig:Mohanty}
\end{figure}
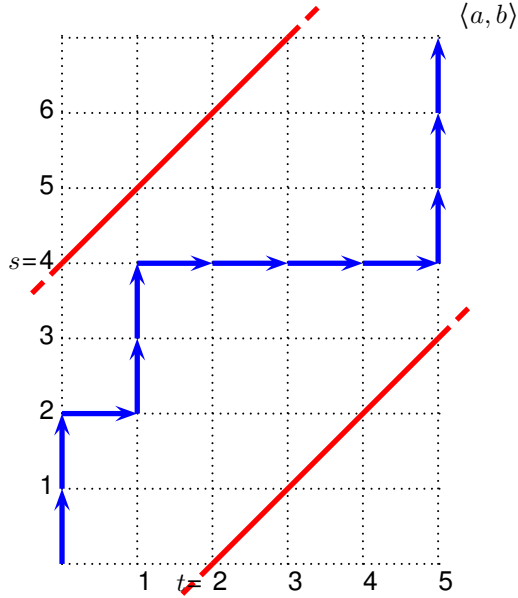

One can derive the enumeration of Theorem~\ref{thm:cnt} using a standard result for bounded lattice paths.
Our main theorem (Theorem~\ref{thm:sym}) will then follow as a corollary.

The number $M(a,b,s,t)$ of
``monotonic'' paths from $\langle 0,0\rangle$ to $\langle a,b\rangle$,  taking $a\in\mathbb N$ steps to the east (E,  $\rightarrow$) and $b\in\mathbb N$ steps to the north (N, $\uparrow$), while totally avoiding (not touching or crossing) the
boundaries $y = x + s$  and $y = x - t$ ($s,t\in\integers^+$),
is known by a reflection argument to be
\begin{align}
M(a,b,s,t) & = \sum_{z \in \integers}\left[\binom{a+b}{b+z(s+t)} - \binom{a+b}{b+z(s+t)+t}\right]
\label{eq:mohanty}
\end{align}
as long as $-t < b-a < s$, so the endpoint $\langle a,b\rangle$ is in bounds.
See \citet[p.\@ 6]{FR,Mohanty}.
An analytic enumeration of corridor paths already appeared in \cite{Ellis}, as pointed out in \cite{CM}.
Figure~\ref{fig:Mohanty} displays a path between boundary lines.

We note that  equation  (\ref{eq:mohanty})  as stated also holds for the cases when $b-a$ is equal to $s$ or to $-t$ since then the two binomials cancel each other for different values of $z$
(for $-z$ when $t=a-b$; 
for $-z-1$ when $s=b-a$).
And indeed there are no admissible paths that end on the boundary.
Thus, we need only limit the formula's applicability to $-t \leq b-a \leq s$.

There is a straightforward relationship between these constrained N/E paths $\langle 0,0\rangle \leadsto\langle a,b\rangle$ and those NE/SE paths $\langle 0,i\rangle \leadsto\langle n,\ell\rangle$ that we have set out to study (as illustrated in Figure~\ref{fig:path}):
\begin{align*}
n &= a+b &
\ell-i &= b-a \\
t &= i+1&
s+t &= h+2 
\end{align*}
Plugging the solution
\begin{align*}
a &= \frac{n+i-\ell}{2} &
b &= \frac{n-i+\ell}{2}\\[1ex]
s &=  h-i+1&
t &= i+1
\end{align*}
into (\ref{eq:mohanty}), 
we get the following:
\begin{align}\label{e}
\AH{i}{\ell}nh
&=\sum_{z \in \integers} \left[\binom{n}{\frac{n-i+\ell}{2}+z(h+2)} - \binom{n}{\frac{n-i+\ell}{2}+z(h+2)+i+1}\right]
\end{align}
as long as $i,\ell\leq h$.
For those $\ell$ for which $\frac{n-i+\ell}{2}$ is not a whole number, the binomial coefficients are taken to be 0, adopting the (nonstandard) convention that $\binom{n}{m}=0$ whenever $m\notin\mathbb N$ \citep[as, e.g., in][p.\@ 75]{Feller}.
See \cite{CM}.

Letting $\ell$ move along the window from $k-j$ to $k+j+1$, we get from  (\ref{e}) that
\begin{align*}
\FFF{i}{k}{j}{n}h
&=\sum_{\substack{\ell=k-j\\0\leq\ell\leq h}}^{k+j+1}
\sum_{z \in \integers}
 \left[\binom{n}{\frac{n-i+\ell}{2}+z(h+2)} - \binom{n}{\frac{n-i+\ell}{2}+z(h+2)+i+1}\right]
\end{align*}
Skipping over the impossible odd or even values of $\ell$ (for which the lower indices of the binomial coefficients are fractional) and shifting summation index, this is:
\begin{align}\label{eq:cnt}
\sum_{z \in \integers}
\sum_{\substack{s=0\\0\leq 2s+p \leq h}}^{j}
\Bigg[ \binom{n}{r+z(h+2)+s}  - \binom{n}{r+z(h+2)+i+s+1} \Bigg]
\end{align}
where $q=n+k-i-j$, $r=\ceil{q/2}$, $p=k-j+q\bmod 2$. 
The different values $s$ can take on correspond to the (at most) $j+1$ possible endpoints in the window $\br{k}{j}$.
The sum for $z$ can be restricted to go from $v=-\floor{(r+i+j+1)/(h+2)}$ to $u=\floor{(n-r)/(h+2)}$
(or to laxer bounds, of course).
We have, thus, arrived at the stated formula of Theorem~\ref{thm:cnt}.

Hereon, consider  only the cases considered in Theorem~\ref{thm:sym}, which guarantee that 
$ k-j\geq -1$ and that $k+j\leq h+1$, so $s$ may run from $0$ to $j$ without exception,
being that $-1\leq 2s+p \leq h+1$ (all that is actually needed) necessarily holds.

When $i=0$,
the inner sum in (\ref{eq:cnt}) simplifies,
leaving only
\begin{align*}
\sum_{z=v}^{u}
\Bigg[ \binom{n}{\ceil{\frac{n+k-j}2}+z(h+2)}  - \binom{n}{\ceil{\frac{n+k+j}2}+z(h+2)+1} \Bigg]
\end{align*}
When also $k=j=h \div 2$, this formula counts Dyck prefixes of bounded height:
\begin{align}\label{eq:dyck}
\sum_{z=v}^{u}
\Bigg[ \binom{n}{\ceil{n/2}+z(h+2)}  - \binom{n}{\ceil{n/2}+\floor{h/2}+z(h+2)+1} \Bigg]
\end{align}
recovering (an equivalent of) the formula in \cite[eq.\@ 1.1]{Cigler}, derived by inclusion/exclusion.
It likewise counts bounded Grand Dyck paths.

When $h$ is $n$ or more, the boundary has no impact on the allowed paths.
Only $z=0$ and the first binomial in (\ref{eq:dyck}) contribute to the count,
yielding simply the central binomial coefficient
\begin{align*}
\binom{n}{\ceil{n/2}} &= \binom{n}{\floor{n/2}}
\end{align*}
which enumerates unrestricted Dyck prefixes, as is well known.

When $k=j=0$ and $n=2m$, these are even-length Dyck paths of semi-length $m$
and bounded height $h$:
\begin{align*}
\sum_{z=-\floor{m/(h+2)}}^{\floor{m/(h+2)}}
\Bigg[ \binom{2m}{m+z(h+2)} - \binom{2m}{m+z(h+2)+1} \Bigg]
\end{align*}
an enumeration due to Howard Grossman in his ``Fun with lattice points'' series \cite[]{fun}.%
\footnote{See sequences 
\seqnum{A011782},
\seqnum{A001519},
\seqnum{A124302},
\seqnum{A080937},
\seqnum{A024175},
\seqnum{A080938},
\seqnum{A033191},
\seqnum{A211216}
 in \citet{OEIS} for $h=2,3,4,5,6,7,8,9$,
 and \seqnum{A080936} for a table for all $h$.}
The same formula counts Dyck prefixes of odd length $2m-1$ ending at $y=1$.

Reversing the order of  summation for the second binomial in (\ref{eq:cnt}), replacing $s$ with $j-s$, we have
\begin{align*}
\sum_{z \in \integers} 
 \Bigg[ 
  \sum_{s=0}^j 
\binom{n}{r+z(h+2)+s}  -
 \sum_{s=0}^j 
 \binom{n}{r+z(h+2)+i+j-s+1}
\Bigg]
\end{align*}
When $j>i$, the inner sums overlap for all $s>i$ and cancel each other,
so one needs only sum for $s$ until $\min\{i,j\}$:
\begin{align*}
\sum_{z \in \integers} 
  \sum_{s=0}^{\min\{i,j\}}
   \Bigg[ 
\binom{n}{r+z(h+2)+s}  -
 \binom{n}{r+z(h+2)+i+j-s+1}
\Bigg]
\end{align*}
As this  is symmetric in $i$ and $j$,
  Theorem~\ref{thm:sym} is vindicated.

\section{Bijective Proof}\label{sec:bi}

A bijection can be inferred from the inductive proof of Section~\ref{sec:ind}
for the equivalence of the enumerations:
\begin{align*}
\FFF{i}{k}{j}{n}{} &= \FFF{j}{k}{i}{n}{} \tag{$\ast$}
\end{align*}
We use a novel representation for paths, which simplifies matters greatly.

Draw a line $y=k+1/2$; we'll call it the \emph{center of attraction}.
Each step starting out toward that line is labeled \texttt{T};
each heading away is labeled \texttt{A}.
From any given point, exactly one outgoing step (\NE\@ or \SE) will be  \texttt{T}
and one \texttt{A}.
Usually, going backward along a  \texttt{T} step is like an  \texttt{A} step,
except when crossing the center line, where it is  \texttt{T} both ways.
Going back  along an \texttt{A} step is always  \texttt{T}.
We call this the \texttt{TA} representation of a lattice path (relative to $k$).
See Figure~\ref{fig:ta}.

Suppose the ordinate of a point along the path is in the window $\br{k}{j}=[k-j:k+j+1]$.
If we take an \texttt{A} step from there, then the next point is in the wider window  $\br{k}{j+1}=[k-j-1:k+j+2]$; 
so $j$ has been incremented.
Conversely, a \texttt{T} step brings it into the narrower range  $\br{k}{j-1}=[k-j+1:k+j]$, with decremented $j$,
unless $j$ is already 0, in which case it stays the same.

If we take this point of view and go through the cases of the inductive proof, 
we find that the correspondence simply reverses the order of steps,
either moving the last step to the beginning or vice versa.
When $i=j$, there is no need to do anything, since the two sides of the equivalence ($\ast$) are identical.
We are led to the following  bijection
between a path $P$ starting at $y=i$ and ending in the range $\br{k}{j}$
and its counterpart path $P^*$ starting at $y=j$ and ending in $\br{k}{i}$:
\begin{itemize}
\item[($=$)]
If $i = j$, then $P^*=P$.
\item[($<$)]
If $i < j$, follow the path from the start at level $i$ until it reaches $j$, if ever.
At that point, we have $P=QR$, where $y=j$ is first reached at the end of prefix $Q$.
Then $P^* = R \vec Q$, where $\vec Q$ is the reverse sequence of $Q$ in its \texttt{TA} representation.
If level $j$ is never attained, then $R$ is empty, and $P^* = \vec P$.
\item[($>$)]
If $i > j$, follow the path from the end backward, starting with a target window of size $j$, 
moving leftwards until it grows to be $i$, if ever.
A \texttt{T} step taken backward enlarges the window, while \texttt{A} shrinks it.
If $R$ is the shortest suffix such that the window size is $i$ at its onset,
so that we have $P=QR$, then
we let  $P^* = \vec R Q$.
If the window never attains size $i$, then $P^* = \vec P$.
\end{itemize}

\begin{figure}[t!]
\centering
\begin{pspicture}(0,-0.5)(12,6.5)
\psgrid[subgriddiv=1,griddots=10](0,0)(0,0)(12,5)
\uput{10pt}[150](-0.1,4.8){$h=$}
\uput{10pt}[150](-0.1,2.8){$k=$}
\uput{10pt}[150](-0.1,0.8){$i=$}
\uput{10pt}[150](-0.1,1.8){$j=$}
\uput{10pt}[150](12.6,4.8){$\ell$}
\uput{15pt}[0](11.7,3.5){$\left.\rule{0pt}{28pt}\right\} 0$}
\uput{15pt}[0](12.2,3.5){$\purple\left.\rule{0pt}{55pt}\right\} 1$}
\uput{15pt}[0](12.7,3.5){$\blue\left.\rule{0pt}{85pt}\right\} 2$}
\psset{linecolor=red,linewidth=2pt,linestyle=solid}
\pcline{-}(0,0)(12,0)
\pcline{-}(0,5)(12,5)
\pcline[linestyle=dashed,linecolor=Green3,linewidth=2pt]{-}(0,3.5)(12,3.5)
\psset{linecolor=blue,arrowsize=6pt,linestyle=solid,linewidth=2pt,labelsep=0pt}
\pcline{->}(0,1)(1,0)\naput{\small\sf\blue A}
\pcline{->}(1,0)(2,1)\naput{\small\sf\blue T}
\pcline{->}(2,1)(3,0)\naput{\small\sf\blue A}
\pcline{->}(3,0)(4,1)\naput{\small\sf\blue T}
\pcline{->}(4,1)(5,0)\naput{\small\sf\blue A}
\pcline{->}(5,0)(6,1)\naput{\small\sf\blue T}
\pcline{->}(6,1)(7,2)\nbput{\small\sf\blue T}
\psset{linecolor=blue,arrowsize=6pt,linestyle=dotted,linewidth=2.5pt,labelsep=0pt}
\pcline{->}(7,2)(8,3)\naput[labelsep=-11pt]{\raisebox{12pt}{\small\sf\blue T}}
\pcline{->}(8,3)(9,4)\naput{\small\sf\blue T}
\pcline{->}(9,4)(10,5)\naput{\small\sf\blue A}
\pcline{->}(10,5)(11,4)\naput{\small\sf\blue T}
\pcline{->}(11,4)(12,5)\naput{\small\sf\blue A}
\psset{linecolor=red!50!blue,arrowsize=6pt,linestyle=dotted,linewidth=2.5pt}
\pcline{->}(0,2)(1,3)\nbput{\small\sf\purple T}
\pcline{->}(1,3)(2,4)\nbput{\small\sf\purple T}
\pcline{->}(2,4)(3,5)\naput{\small\sf\purple A}
\pcline{->}(3,5)(4,4)\naput{\small\sf\purple T}
\pcline{->}(4,4)(5,5)\naput{\small\sf\purple A}
\psset{linecolor=red!50!blue,arrowsize=6pt,linestyle=solid,linewidth=2pt}
\pcline{->}(5,5)(6,4)\naput{\small\sf\purple T}    
\pcline{->}(6,4)(7,3)\naput{\small\sf\purple T}
\pcline{->}(7,3)(8,2)\nbput[labelsep=-10pt]{\raisebox{-19pt}{\small\sf\purple A}}
\pcline{->}(8,2)(9,3)\naput{\small\sf\purple T}
\pcline{->}(9,3)(10,2)\nbput{\small\sf\purple A}
\pcline{->}(10,2)(11,3)\nbput{\small\sf\purple T}
\pcline{->}(11,3)(12,2)\nbput{\small\sf\purple A}
\psset{linecolor=brown,arrowsize=3pt,linestyle=solid,linewidth=1pt}
\pcline{|-|}(12,1.7)(12,4.3)
\pcline{|-|}(11,2.7)(11,3.3)
\pcline{|-|}(10,1.7)(10,4.3)
\pcline{|-|}(9,2.7)(9,3.3)
\pcline{|-|}(8,1.7)(8,4.3)
\pcline{|-|}(7,2.7)(7,3.3)
\pcline{|-|}(6,1.7)(6,4.3)
\pcline{|-|}(5,0.7)(5,5.3)
\end{pspicture}
\medskip
\caption{The 12-step (blue) path $P$ from level $i=1$ to  level $\ell=5$
belongs to the class enumerated by $\FFF{1}{3}{2}{12}{}$.
Since $k=3$,
steps are labeled \textsf{T} when they start out toward the ``attractor'' $y=3.5$ (the dashed {\color{Green4}green} line), and \textsf{A} when they head away in the opposite direction.
Accordingly, path $P$ is labeled \textsf{\blue ATATATTTTATA}.
Windows of sizes 0, 1, and 2 around the green attractor are shown to the right of the grid.
The target window size is $j=2$, so we are in the ($<$) case of the bijection.
After seven (solid blue) steps \textsf{\blue ATATATT}, the path touches $y=2$, 
so the remaining 5 (dotted) steps, 
\textsf{\blue TTATA}, are copied as is and placed with $\langle 0,2\rangle$ as their initial point (dotted purple),
followed by the first seven in reverse (solid purple), that is, \textsf{\purple TTATATA}, to obtain the corresponding path $P^*$.
The result is \textsf{\purple TTATATTATATA}, one of those counted by $\FFF{2}{3}{1}{12}{}$,
which all start from $j=2$ and end in $[2\dotdot 5]$.
The lower path $P$ minus its head, viz.\@ $\blue \A{0}{5}{11}$, corresponds to the upper (purple) path $P^*$ minus its tail, $\purple\A{2}{3}{11}$, which ends in a size 0 window.
The initial \textsf{\blue A} step serves to extend that 11-step (purple) path with a downwards \textsf{\purple A} step that enlarges the window but remains in bounds.
The counterpart of this (purple) path is again the original (in blue), and is obtained
by proceeding from the end toward the beginning until the window size becomes 2,
per case ($>$).
The actual windows at each point along the reversed (solid purple) segment are shown (as brown bars),
going from size 2 at $x=5$, through 1, 0, 1, 0, 1, 0,  ending with size 1.
On account of the unusual encoding, the reversed (solid) path segments do not actually resemble each other visually.
See Section~\ref{sec:bi} for details.}\label{fig:ta}
\end{figure}

The path in Figure \ref{fig:path} is its own counterpart, as this is an instance of case ($=$) with $i=j=1$.
For a worked-out nontrivial example, see Figure \ref{fig:ta}.

Consider case ($<$), when $i<j$.
A path in $\FF{i}{k}{j}{n}$ is composed of two parts: the $Q$ part starting from $y=i$ and never crossing $y=j$;
and the $R$ part, which~-- when nonempty~-- goes from $y=j$ until $y=\ell\in\br{k}{j}$.
Considering that $i<j\leq k+1$, all \texttt{T} steps in $Q$ are upward, and all \texttt{A}'s are downward.
So the net difference in ordinate value at the start and end of $Q$ is at most $j-i$,
and likewise the number of \texttt{T}'s in $Q$ minus the number of \texttt{A}'s is at most $j-i$.
The corresponding segment ${\vec Q}$ in the path $P^* = R {\vec Q}$ starts either at $\ell$ in the nonempty case or else at $j$.
In either event, the starting point of ${\vec Q}$ is within $\br{k}{j}$.
In the nonempty case, because $R$ ends at $\ell$;
in the empty case, on account of the fact that $j$ is higher than the end of $P$, which itself is in $\br{k}{j}$.
Each \texttt{T} step in ${\vec Q}$ moves the path toward the center of attraction and shrinks the window; 
each \texttt{A} step enlarges it.
Hence, the ending point is within $\br{k}{j-(j-i)}=\br{k}{i}$, as required.

It is also not hard to verify that the transpositions involved keep the path within the bounded corridor, given that the original path
satisfies $0\leq i,j\leq\min\{k+1,h-k\}$.
The points along the $Q$ segment all lie in the range $[0\dotdot j]$,
reaching $j$ only at the end~-- if ever.
Thus, the corresponding path stays within the window $\brk{k\pm j}$,
and can only be at the edge of that window initially.
Just as each \texttt{T} step increases the vertical position (since $k$ lies above) along the given path
but never exceeds $j$,
so too the \texttt{T} step in ${\vec Q}$ shrinks the window from its initial maximal size $j$,
but never below 0.
Likewise, an \texttt{A} step decreases the vertical position down to $0$,
while the \texttt{A} step in the corresponding path enlarges the window, but only up to size $j$.
Those remaining steps $R$ that are copied as is clearly remain in bounds. 
Figure~\ref{fig:ta} portrays how the window size changes along ${\vec Q}$.

The symmetric case, with $j<i$, is perfectly analogous, starting instead from the paths $\FF{j}{k}{i}{n}$ counted on the right side of ($\ast$).

Thus we have attained a bijective proof of Theorem~\ref{thm:sym}.
We note that the inductive proof allows for alternate bijections depending on the preferred order in which the different cases are to be considered.

\section{Historical Discussion}

Theorem~\ref{thm:sym}, our main result, is a significant generalization of the equality given by Johann Cigler \citeyearpar{Cigler}, namely,
\begin{align}
\label{c} \AH{0}{[0\dotdot h]}{n}h &~~=~~ \AH{\hslash}{[\hslash\dotdot \hslash+1]}{n}h
\end{align}
for all heights $h$, where
$\hslash=h\div 2$ for short.
Paths (counted by) $\A{\hslash}{[\hslash\dotdot \hslash+1]}{n} $ start in the middle of the swath and end either in the middle~-- when the number of steps is even, or just above~-- when odd.
As noted earlier, these are called ``Grand Dyck'' paths.
Dyck path prefixes
$ \A{0}{[0\dotdot h]}{n}$ start at the bottom and end anywhere within the swath.
Cigler's (\ref{c}) asserts the equality of cardinality of these two sets of paths.
As such, it
is a particular instance of our more general result (\ref{eq:sym}) with $i=0$ and $j=k=\hslash$.
Phrased in our notation, Cigler demonstrated:
\begin{align*}
\FFF{0}{\hslash}{\hslash}{n}h & ~~=~~ \FFF{\hslash}{\hslash}0{n}h
\end{align*}

Cigler solicited alternative proofs of his result.
More specifically, he asked in \citet{Cigler15} for a \emph{bijective} proof of the height $h=3$ case,
\begin{align}\nonumber
 \GH{0}{0\dotdot 3}{n}3 &~~=~~ \AH{1}{[1\dotdot 2]}{n}3
\end{align}
which gives rise to the Fibonacci numbers.
The wished-for bijective solution to this very particular case was discovered shortly thereafter by Thomas Prellberg \cite[Answer]{Cigler15}, followed by another due to Helmut Prodinger \citeyearpar{prodinger2016height}.
Most recently, Nancy Gu and Prodinger \citeyearpar{GP} constructed a bijection for Cigler's full case (\ref{c}) by extending the idea in \citet{prodinger2016height}.

When there are no upper and lower bounds on paths, there are long-standing, well-known bijections  between
Grand Dyck paths and Dyck path prefixes,
\begin{align*}
 \FFF{0}{k}{k}{n}\infty & ~~=~~ \FFF{k}{k}{0}{n}\infty 
\end{align*}
for $k\geq n/2$,
as mentioned \citep{Feller,Greene}.
Endre Cs\'aki and Sri Gopal Mohanty \citeyearpar[Theorem 3.1]{CM} , in fact, already proved 
the following generalization to corridors by means of an inductively constructed bijection:
\begin{align*}
\FFF{0}{k}{k}{n}h & ~~=~~ \FFF{k}{k}{0}{n}h
\end{align*}

Our bijection in the previous section supplies an alternative proof of Cigler's instance (\ref{c}).
In that special case,
the bijection amounts to simply reversing the order of steps in the \texttt{TA} representation.
This works as is for even $n$ in the Grand Dyck ($i=k=\hslash$ and $j=0$)
to Dyck-prefix ($i=0$ and $j=k=\hslash$) case of Cigler,
as this is the ($i>j$) case of the bijection and the window never grows too big 
to continue all way to the beginning.
(It may get to be $\hslash$; the maximum excess of \texttt{T} moves over \texttt{A} moves, but no larger.)
Unfortunately, simple reversal doesn't do the trick when $n$ is odd and $h$ is even 
because proceeding only backward can lead 
to a window wider than $\hslash=i$.
(For example, \texttt{TAT} is a valid path in $\G{0}{0\dotdot 3}{3}$, but \texttt{TAT} is not in $\G{1}{1\dotdot 2}{3}$. 
The bijection of the previous section reverses only part and yields the valid \texttt{ATT} instead.)
For the odd $n$ case,
it is possible to modify the bijection by first reversing the Grand Dyck path left to right
(so it ends on $y=\hslash$ but begins at $i=\hslash+1$)
before converting to the \texttt{TA} representation and reversing.
This now covers all cases of (\ref{c}).
The second bijection also works for even $h$ and even $n$.
For odd $h$, regardless of the parity of $n$, the first bijection actually succeeds for all $i,j$ meeting the requirements of the theorem.
So, when $n$ and $h$ have the same parity, both bijections work.
In the more general cases, when $k\neq\hslash$, neither applies, and we resort to 
the slightly more complicated bijection of the previous section, wherein only part of the \texttt{TA} path is reversed.

We began our investigation seeking a bijective proof of (\ref{c}).
The simple bijection employing the \texttt{TA} path encoding didn't work in all cases.  
This led us to a sequence of generalizations,
commencing from Cigler's (\ref{c}):
\begin{align*}
\G{0}{0\dotdot h}{n} &= \G{\hslash}{\hslash\dotdot \hslash+1}{n}\\
 \G{i}{0\dotdot h}{n} &= \G{\hslash}{\hslash-i\dotdot \hslash +i+1}{n}\\
 \G{i}{\hslash-j\dotdot \hslash+j+1}{n} &= \G{j}{\hslash-i\dotdot \hslash+i+1}{n}\\
 \G{i}{k-j\dotdot k+j+1}{n} &= \G{j}{k-i\dotdot k+i+1}{n}
\end{align*}
First we let $i$ be anywhere (not just $0$), then we let $j$ be any size (not just $\hslash$), and finally allowed it to be centered at any $k$ (not just $\hslash$).
Concurrently,
we programmed various enumerations and potential bijections to lend support to~-- or refute~--  conjectures as they arose.
Casting the equivalence in a fashion that highlights its symmetry also contributed to finding the 
generalizations and proofs. 

All the above variants share the basic idea that, as the starting point of one set of paths moves from the edge of the corridor toward the middle, the target range of the corresponding equinumerous set of paths grows wider and wider.
This behavior is what suggested the \texttt{TA} encoding in the first place.

\acknowledgments
I gratefully thank Johann Cigler for  encouragement, references, and advice, Christian Rinderknecht for first bringing Cigler's interesting challenge  to my attention, and a referee for her critique.

\bibliographystyle{abbrvnat}
\bibliography{Paths}

\clearpage

\begin{table}[p]
\centering\tiny\setlength\tabcolsep{2pt}\renewcommand\arraystretch{1.5}
\begin{tabular}{|r||d|d|d|d|d|d|d|d|d|d|d|d|d|d|d|d|d||l||c|}
$i$ & \multicolumn{1}{c|}{$n=0$} & \multicolumn{1}{c|}{1} & \multicolumn{1}{c|}{2} & \multicolumn{1}{c|}{3} & \multicolumn{1}{c|}{4} & \multicolumn{1}{c|}{5} & \multicolumn{1}{c|}{6} & \multicolumn{1}{c|}{7} & \multicolumn{1}{c|}{8} & \multicolumn{1}{c|}{9} & \multicolumn{1}{c|}{10} & \multicolumn{1}{c|}{11} & \multicolumn{1}{c|}{12} & \multicolumn{1}{c|}{13} & \multicolumn{1}{c|}{14} & \multicolumn{1}{c|}{15} & \multicolumn{1}{c||}{16} & $\ell$ & OEIS\\\hline\hline
4 & 1 & &1 & &2 & &5 & &14 & &41 & &122 & &365 & &1094 &4  & \seqnum{A007051} \\\hline
 &  &1 & &2 & &5 & &14 & &41 & &122 & &365 & &1094 & &3 & \seqnum{A007051}\\\hline
 &  & &1 & &3 & &9 & &27 & &81 & &243 & &729 & &2187 &2 & \seqnum{A000244}\\\hline
 &  & & &1 & &4 & &13 & &40 & &121 & &364 & &1093 & &1 & \seqnum{A003462}\\\hline
 &  & & & &1 & &4 & &13 & &40 & &121 & &364 & &1093 &0 & \seqnum{A003462}\\\hline
\hline
 &  &1 & &2 & &5 & &14 & &41 & &122 & &365 & &1094 & &4 & \seqnum{A007051}\\\hline
3 & 1 & &2 & &5 & &14 & &41 & &122 & &365 & &1094 & &3281 &3 & \seqnum{A007051}\\\hline
 &  &1 & &3 & &9 & &27 & &81 & &243 & &729 & &2187 & &2 & \seqnum{A000244}\\\hline
 &  & &1 & &4 & &13 & &40 & &121 & &364 & &1093 & &3280 &1 & \seqnum{A003462}\\\hline
 &  & & &1 & &4 & &13 & &40 & &121 & &364 & &1093 & &0 & \seqnum{A003462}\\\hline
\hline
 &  & &1 & &3 & &9 & &27 & &81 & &243 & &729 & &2187 &4 & \seqnum{A000244}\\\hline
 &  &1 & &3 & &9 & &27 & &81 & &243 & &729 & &2187 & &3 & \seqnum{A000244}\\\hline
2 & 1 & &2 & &6 & &18 & &54 & &162 & &486 & &1458 & &4374 &2 & \seqnum{A025192}\\\hline
 &  &1 & &3 & &9 & &27 & &81 & &243 & &729 & &2187 & &1 & \seqnum{A000244}\\\hline
 &  & &1 & &3 & &9 & &27 & &81 & &243 & &729 & &2187 &0 & \seqnum{A000244}\\\hline
\hline
 &  & & &1 & & \multicolumn{1}{r|}{\bf\blue 4} & &13 & &40 & &121 & &364 & &1093 & &4 & \seqnum{A003462}\\\hline
 &  & & \multicolumn{1}{r|}{\bf\blue 1} & &  \multicolumn{1}{r|}{\bf\blue 4} & & \multicolumn{1}{r|}{\bf\blue 13} & &40 & & 121 & &\multicolumn{1}{r|}{\bf\blue 364} & &1093 & &3280 &3& \seqnum{A003462} \\\hline
 &  &  \multicolumn{1}{r|}{\bf\blue 1} & & \multicolumn{1}{r|}{\bf\blue 3} & &9 & & \multicolumn{1}{r|}{\bf\blue 27}& & 81 & &\multicolumn{1}{r|}{\bf\blue 243} & &729 & &2187 & &2 & \seqnum{A000244}\\\hline
1 & \multicolumn{1}{r|}{\bf\blue 1} & &2 & &5 & &14 & & \multicolumn{1}{r|}{\bf\blue 41} & &\multicolumn{1}{r|}{\bf\blue 122} & &365 & &1094 & &3281 &1& \seqnum{A007051} \\\hline
 &  &1 & &2 & &5 & &14 & &\multicolumn{1}{r|}{\bf\blue 41} & &122 & &365 & &1094 & &0 & \seqnum{A007051}\\\hline
\hline
 &  & & & &1 & &4 & &13 & &40 & &121 & &364 & &1093 &4 & \seqnum{A003462}\\\hline
 &  & & &1 & &4 & &13 & &40 & &121 & &364 & &1093 & &3& \seqnum{A003462} \\\hline
 &  & &1 & &3 & &9 & &27 & &81 & &243 & &729 & &2187 &2 & \seqnum{A000244}\\\hline
 &  &1 & &2 & &5 & &14 & &41 & &122 & &365 & &1094 & &1 & \seqnum{A007051}\\\hline
0 & 1 & &1 & &2 & &5 & &14 & &41 & &122 & &365 & &1094 &0 & \seqnum{A007051}\\\hline
\end{tabular}
\vspace*{1em}
\caption{The number of paths $\AH{i}{\ell}{n}{4}$ for $i,\ell\in[0\dotdot 4]$, $n\in[0\dotdot 16]$,
is found in the subtable labeled $i$ (in the first column), column $n$, and row $\ell$ (penultimate column) of that subtable.
For example, $\A{2}{2}{16}=\A{0}{[0\dotdot 4]}{16}=4374$ 
and 
 $\A{3}{3}{16}=\A{4}{[2\dotdot 4]}{16}=3281$.
The former corresponds to $\FFF{2}{2}{0}{}{}=\FFF{0}{2}{2}{}{}$ per Theorem \ref{thm:sym};
the latter follows from  $\FFF{3}{3}{0}{}{}=\FFF{1}{0}{0}{}{}$ by Lemma \ref{lem} and $\FFF{1}{0}{0}{}{}=\FFF{4}{3}{1}{}{}$ by Corollary \ref{cor:2}.
As for a bishop on a chessboard, half the squares are unreachable from any given starting point;
the few squares that require backward steps are likewise inaccessible.
The particular path of Figures~\ref{fig:path}--\ref{fig:Mohanty} is included in those counted by the numbers highlighted in \textbf{\blue blue boldface}.
Sloane (OEIS) numbers of the row sequences are provided in the last column.}\label{tab}
\end{table}

\begin{table}[p]
\centering\tiny\setlength\tabcolsep{2pt}\renewcommand\arraystretch{1.5}
\begin{tabular}{|r||d|d|d|d|d|d|d|d|d|d|d|d|d|d|d|d|d||l||c|}
$i$ & \multicolumn{1}{c|}{$n=0$} & \multicolumn{1}{c|}{1} & \multicolumn{1}{c|}{2} & \multicolumn{1}{c|}{3} & \multicolumn{1}{c|}{4} & \multicolumn{1}{c|}{5} & \multicolumn{1}{c|}{6} & \multicolumn{1}{c|}{7} & \multicolumn{1}{c|}{8} & \multicolumn{1}{c|}{9} & \multicolumn{1}{c|}{10} & \multicolumn{1}{c|}{11} & \multicolumn{1}{c|}{12} & \multicolumn{1}{c|}{13} & \multicolumn{1}{c|}{14} & \multicolumn{1}{c|}{15} & \multicolumn{1}{c||}{16} & $\ell$ & OEIS\\\hline\hline
5 & 1 & &1 & &2 & &5 & &14 & &42 & &131 & &417 & &1341 &5  & \seqnum{A080937}\\\hline
 &  &1 & &2 & &5 & &14 & &42 & &131 & &417 & &1341 & &4 & \seqnum{A080937}\\\hline
 &  & &1 & &3 & &9 & &28 & &89 & &286 & &924 & &2993 &3 & \seqnum{A094790}\\\hline
 &  & & &1 & &4 & &14 & &47 & &155 & &507 & &1652 & &2 & \seqnum{a094789}\\\hline
 &  & & & &1 & &5 & &19 & &66 & &221 & &728 & &2380 &1 & \seqnum{A005021}\\\hline
 &  & & & & &1 & &5 & &19 & &66 & &221 & &728 & &0 & \seqnum{A005021}\\\hline
\hline
 &  &1 & &2 & &5 & &14 & &42 & &131 & &417 & &1341 & &5 & \seqnum{A080937}\\\hline
4 & 1 & &2 & &5 & &14 & &42 & &131 & &417 & &1341 & &4334 &4 & \seqnum{A080937}\\\hline
 &  &1 & &3 & &9 & &28 & &89 & &286 & &924 & &2993 & &3 & \seqnum{A094790}\\\hline
 &  & &1 & &4 & &14 & &47 & &155 & &507 & &1652 & &5373 &2 & \seqnum{a094789}\\\hline
 &  & & &1 & &5 & &19 & &66 & &221 & &728 & &2380 & &1 & \seqnum{A005021}\\\hline
 &  & & & &1 & &5 & &19 & &66 & &221 & &728 & &2380 &0 & \seqnum{A005021}\\\hline
\hline
 &  & &1 & &3 & &9 & &28 & &89 & &286 & &924 & &2993 &5 & \seqnum{A094790}\\\hline
 &  &1 & &3 & &9 & &28 & &89 & &286 & &924 & &2993 & &4 & \seqnum{A094790}\\\hline
3 & 1 & &2 & &6 & &19 & &61 & &197 & &638 & &2069 & &6714 &3 & \seqnum{A052975}\\\hline
 &  &1 & &3 & &10 & &33 & &108 & &352 & &1145 & &3721 & &2 & \seqnum{A060557}\\\hline
 &  & &1 & &4 & &14 & &47 & &155 & &507 & &1652 & &5373 &1 & \seqnum{a094789}\\\hline
 &  & & &1 & &4 & &14 & &47 & &155 & &507 & &1652 & &0 & \seqnum{a094789}\\\hline
\hline
 &  & & &1 & &4 & &14 & &47 & &155 & &507 & &1652 & &5 & \seqnum{a094789}\\\hline
 &  & &1 & &4 & &14 & &47 & &155 & &507 & &1652 & &5373 &4 & \seqnum{a094789}\\\hline
 &  &1 & &3 & &10 & &33 & &108 & &352 & &1145 & &3721 & &3 & \seqnum{A060557}\\\hline
2 & 1 & &2 & &6 & &19 & &61 & &197 & &638 & &2069 & &6714 &2 & \seqnum{A052975}\\\hline
 &  &1 & &3 & &9 & &28 & &89 & &286 & &924 & &2993 & &1 & \seqnum{A094790}\\\hline
 &  & &1 & &3 & &9 & &28 & &89 & &286 & &924 & &2993 &0 & \seqnum{A094790}\\\hline
\hline
 &  & & & &1 & &5 & &19 & &66 & &221 & &728 & &2380 &5 & \seqnum{A005021}\\\hline
 &  & & &1 & &5 & &19 & &66 & &221 & &728 & &2380 & &4 & \seqnum{A005021}\\\hline
 &  & &1 & &4 & &14 & &47 & &155 & &507 & &1652 & &5373 &3 & \seqnum{a094789}\\\hline
 &  &1 & &3 & &9 & &28 & &89 & &286 & &924 & &2993 & &2 & \seqnum{A094790}\\\hline
1 & 1 & &2 & &5 & &14 & &42 & &131 & &417 & &1341 & &4334 &1 & \seqnum{A080937}\\\hline
 &  &1 & &2 & &5 & &14 & &42 & &131 & &417 & &1341 & &0 & \seqnum{A080937}\\\hline
\hline
 &  & & & & &1 & &5 & &19 & &66 & &221 & &728 & &5 & \seqnum{A005021}\\\hline
 &  & & & &1 & &5 & &19 & &66 & &221 & &728 & &2380 &4 & \seqnum{A005021}\\\hline
 &  & & &1 & &4 & &14 & &47 & &155 & &507 & &1652 & &3 & \seqnum{a094789}\\\hline
 &  & &1 & &3 & &9 & &28 & &89 & &286 & &924 & &2993 &2 & \seqnum{A094790}\\\hline
 &  &1 & &2 & &5 & &14 & &42 & &131 & &417 & &1341 & &1 & \seqnum{A080937}\\\hline
0 & 1 & &1 & &2 & &5 & &14 & &42 & &131 & &417 & &1341 &0 & \seqnum{A080937}\\\hline
\end{tabular}
\vspace*{1em}
\caption{Paths $\AH{i}{\ell}{n}{5}$ constrained to height 5, $n\in[0\dotdot 16]$.
Sloane (OEIS) numbers of the row sequences are provided in the last column.}\label{tab:5}
\end{table}

\end{document}